\documentclass{amsart}

\usepackage{xypic}
\input xy
\xyoption{all}
\usepackage{epsfig}
\usepackage{amsthm}
\usepackage{amssymb}
\usepackage{amsmath}
\usepackage{amscd}

\newcommand{\bg}{\begin{equation}}
\newcommand{\ed}{\end{equation}}
\newcommand{\bga}{\begin{eqnarray}}
\newcommand{\eda}{\end{eqnarray}}
\newcommand{\pf}{\textbf{Proof:\ }}
\newcommand{\lpf}{\textbf{Proof of Lemma:\ }}

\def\cbdu{\par{\raggedleft$\Box$\par}}

\newtheorem {Theorem}  {Theorem}

\numberwithin{Theorem}{section}

\newtheorem {Lemma}[Theorem]  {Lemma}

\theoremstyle{definition}
\newtheorem{Definition}[Theorem]{Definition}
\theoremstyle{remark}
\newtheorem{Remark}[Theorem]{Remark}

\newtheorem {Corollary}[Theorem]{Corollary}
\newtheorem{Claim}[Theorem]{Claim}
\newtheorem{Assumption}[Theorem]{Assumption}
\expandafter\chardef\csname pre amssym.def
at\endcsname=\the\catcode`\@ \catcode`\@=11
\def\undefine#1{\let#1\undefined}
\def\newsymbol#1#2#3#4#5{\let\next@\relax
 \ifnum#2=\@ne\let\next@\msafam@\else
 \ifnum#2=\tw@\let\next@\msbfam@\fi\fi
 \mathchardef#1="#3\next@#4#5}
\def\mathhexbox@#1#2#3{\relax
 \ifmmode\mathpalette{}{\m@th\mathchar"#1#2#3}%
 \else\leavevmode\hbox{$\m@th\mathchar"#1#2#3$}\fi}
\def\hexnumber@#1{\ifcase#1 0\or 1\or 2\or 3\or 4\or 5\or 6\or 7\or 8\or
 9\or A\or B\or C\or D\or E\or F\fi}

\font\teneufm=eufm10 \font\seveneufm=eufm7 \font\fiveeufm=eufm5
\newfam\eufmfam
\textfont\eufmfam=\teneufm \scriptfont\eufmfam=\seveneufm
\scriptscriptfont\eufmfam=\fiveeufm

\catcode`\@=\csname pre amssym.def at\endcsname

\newcounter{remark}
\newcommand{\e}{\epsilon}

\newcommand{\R}{\mathbf{R}}
\renewcommand{\div}{\mbox{div}}
\newcommand{\Aa}{{\mathcal A}}
\newcommand{\Bb}{{\mathcal B}}
\newcommand{\Cc}{{\mathcal C}}

\newcommand{\Dd}{{\mathcal D}}

\newcommand{\Ll}{{\mathcal L}}

\newcommand{\Hh}{{\mathcal H}}

\newcommand{\Mm}{{\mathcal M}}

\def  \R   {{\mathbb R}}

\def  \12  {{\frac{1}{2}}}

\def\build#1_#2^#3{\mathrel{\mathop{\kern 0pt#1}\limits_{#2}^{#3}}}

\begin{document}

\title{Regularity of Solutions to the Liquid Crystals Systems in $\mathbb{R}^2$ and $\mathbb{R}^3$}

\author[Mimi Dai]{ Mimi Dai}
\address{Department of Mathematics, UC Santa Cruz, Santa Cruz, CA 95064,USA}
\email{mdai@ucsc.edu}
\author [Jie Qing]{Jie Qing}
\address{Department of Mathematics, UC Santa Cruz, Santa Cruz, CA 95064,USA}
\email{qing@ucsc.edu}
\author[Maria Schonbek] {Maria Schonbek}
\address{Department of MAthematics, UC Santa Cruz, Santa Cruz, CA 95064, USA}
\email{schonbek@ucsc.edu}


\thanks{The work of Mimi Dai was partially supported by NSF Grant
DMS-0700535 and DMS-0600692;}
\thanks{The work of Jie Qing was
partially supported by NSF Grant DMS-0700535;}
\thanks{The work of
M. Schonbek was partially supported by NSF Grant DMS-0600692}

\begin{abstract}
Global existence for weak solutions to systems of nematic liquid crystals, with non-constant fluid density has been established in \cite{Liu} and \cite{JT}.  In this paper we extend the regularity and uniqueness results of Fanghua Lin and Chun Liu in \cite{LL} for the systems of nematic liquid crystals (LCD). In \cite{LL}, the underlying system has constant density. In this paper, the regularity and uniqueness results are established for a density dependent LCD system.
\end{abstract}

\maketitle

\section{Introduction}

The flows of nematic liquid crystals can be treated as slow moving particles where the fluid velocity and the alignment of the particles influence each other. The hydrodynamic theory of liquid crystals  was established by Ericksen \cite{Er0}, \cite{Er1} and Leslie \cite{Le0}, \cite{Le1} in the 1960's.
As Leslie points out in his 1968 paper : ``liquid crystals are states of matter which are capable of flow, and in which the molecular arrangements  give rise to a preferred direction".  In this paper we consider the simplified model for the flow of nematic liquid crystals:

\begin{equation}\begin{split}\label{LCD}
\rho_t+\nabla\cdot (\rho u) =0,\\
(\rho u)_t+\nabla\cdot(\rho u\otimes u)+\nabla p =\triangle u -\nabla\cdot(\nabla d\otimes\nabla d)\\
d_t+u\cdot\nabla d =\triangle d-f(d)\\
\nabla\cdot u =0
\end{split}
\end{equation}
in $\Omega\times(0, T)$, where $\Omega$ is a domain in $R^n$,
$\rho: \Omega\times[0,T]\to\mathbb{R}$ is the fluid density, $
p: \Omega\times[0,T]\to\mathbb{R}$ is the fluid pressure,  $
u: \Omega\times[0,T]\to\mathbb{R}^n$ is the fluid velocity and
$d: \Omega\times[0,T]\to\mathbb{R}^n$ is the director field representing the alignment of the molecules, with $n=2,3$. The force term $\nabla d\otimes\nabla d$ in the equation of the conservation of momentum denotes the $3\times 3$ matrix whose $ij$-th entry is given by $\nabla_i d\cdot\nabla_j d$ for $1\leq i,j\leq 3$. This force $\nabla d\otimes\nabla d$ is the stress tensor of the energy about the director field $d$, where the energy is given by:
$$
\frac 12 \int_{\Omega} |\nabla d|^2 dx + \int_{\Omega}F(d)dx
$$
and
$$
F(d)=\frac{1}{4\eta^2}(|d|^2-1)^2, \quad f(d) = \nabla F(d) = \frac{1}{\eta^2}(|d|^2-1)d.
$$
In fact $F(d)$ is the penalty term of the Ginzburg-Landau approximation of the original free energy of the director field with unit length.\\

There is a vast literature on the hydrodynamic of the liquid crystal system. For background we list  a few names, with no intention to be complete:  \cite{EK},  \cite{HKL}, \cite{Kin}, \cite{LL}, \cite{LL2}, \cite{LL1}, \cite{CC}, \cite{Cal}, \cite{CC}, \cite{CC1}, \cite{CDLL}, \cite {Liu} and \cite{JT}.  Particularly, in \cite{Liu} \cite{JT}, the global weak existence of solutions to the flow of nematic liquid crystals was obtained for  fluids with non constant  density. In light of the regularity results to the pure fluid system established in \cite{AKM} and \cite{LS1} it is natural to ask if  the regularity results in \cite{LL}   can  be extend to prove the regularity of the solutions for flows of nematic liquid crystals with non-constant fluid density.

In this paper we focus on the regularity of solutions to the flow of nematic liquid crystals satisfying the initial conditions:
\bg\label{initden}
\rho(x,0) =\rho_0(x), \ \ \ 0<M_1\leq\rho_0(x)\leq M_2,
\ed
\bg\label{initu}
u(x,0) =u_0(x), \ \ \ \nabla\cdot u_0 =0, \ \ \ u_0|_{\partial\Omega}=0,
\ed
\bg\label{initd}
d(x,0) =d_0(x), \ \ \ |d_0(x)|= 1,
\ed
and the boundary conditions:
\bg\label{bd}
u(\cdot, t)|_{\partial\Omega}=0, \ \ \ d(\cdot, t)|_{\partial\Omega}=d_0|_{\partial\Omega}.
\ed

Existence of  global weak solutions of (\ref{LCD}),  with the above specified data, has been established in \cite{Liu} and \cite{JT}. In fact they have existence even without assuming the positive lower bound $M_1$.
In  \cite{Liu} and \cite{JT}, to derive global weak solution, a viscosity term $\epsilon\Delta \rho$ is added  to regularize  the first equation of the system \eqref{LCD}. This approach had been  suggested in \cite{Lio}. Our proof of regularity uses  energy estimates  introduced by Ladyzhenskaya  on the approximate solutions. The added term in the second equation that results from the regularizing  viscosity  in the first equation  in  \cite{Lio}, contains the gradient of the fluid density. This term   seems to create difficulties when it is used   to establish the Ladyzhenskaya energy estimates for the approximate solution derived by the Galerkin method. Thus in the appendix we  sketch a proof of  existence for the global weak solutions to  system \eqref{LCD} without  the introduction of the  viscosity  term for the density in the equation of the conservation of mass. In our case to obtain a classical solution we need to work with data that is more regular than the data  used in \cite{Liu} and \cite{JT}.\\

We obtain interior regularity with a relatively weak conditions on the initial data. For more regular data, we are able to obtain solutions which are regular up to the boundary.  In the rest of the introduction we briefly  describe  our main results: \\

\noindent \emph{Regularity in 2 dimensions:}
\begin{Theorem}\label{Mthm1}Suppose that $\Omega$ is a smooth bounded domain in $\R^2$. Let $\rho_0$, $u_0$ and $d_0$ satisfy (\ref{initden})-(\ref{bd}).
Suppose that $\rho_0\in C^1(\Omega)$, $u_0\in H^1(\Omega)$ and $d_0\in H^{2}(\Omega)$. Then, the system (\ref{LCD}) has a global classical solution $(\rho,u,d)$, that is, for all $T>0$ and some $\alpha\in(0,1)$
$$
u\in C^{1+\alpha/2,2+\alpha}((0,T)\times\Omega)
$$
$$
\nabla p\in C^{\alpha/2,\alpha}((0,T)\times\Omega)
$$
$$
d\in C^{1+\alpha/2,2+\alpha}((0,T)\times\Omega)
$$
$$
\rho\in C^1((0,T)\times\Omega).
$$
\end{Theorem}

\begin{Corollary}\label{Cor1}
Suppose $(\rho_0,u_0,d_0)$ satisfies the conditions in Theorem \ref{Mthm1}. In addition, $\rho_0\in C^1(\bar\Omega)$,  $u_0\in C^{2+\alpha}(\bar\Omega)$ and $d_0\in C^{2+\alpha}(\bar\Omega)$. Then the solution $(\rho,u,d)$ is regular up to boundary, that is, the conclusions in Theorem \ref{Mthm1} hold in$[0,T]\times\bar\Omega$.
\end{Corollary}

Indeed, more regularity on the initial data $d_0$ implies more regularity on the boundary due to the second condition in (\ref{bd}). Thus, yielding the regularity up to boundary, see \cite{Kr1}.\\

The regularity for the flow of the nematic liquid crystals in dimension 2 turns out to be not too difficult. We first establish the Ladyzhenskaya energy estimate \eqref{phiest} and \eqref{L^2H^2}, similar to that in \cite{LL}. Then we apply the regularity result for transport equations in \cite{AKM} to obtain the H\"{o}lder continuity of the fluid density. Therefore
Theorem \ref{Mthm1} follows from the $L^p$ estimates and H\"older estimates in \cite{LS} and a more or less standard bootstrapping between the three equations in the
system \eqref{LCD}.\\


Provided we have sufficiently small data or we work with sufficiently short time in 3D, we also obtain interior regularity. Given more restriction on data, the regularity can be obtained up to the boundary.\\

\noindent \emph{Regularity  in  3 dimension  with small data or short time:}
\begin{Theorem}\label{Mthm2} Suppose that $\Omega$ is a smooth bounded domain in $\R^3$. Let $\rho_0$, $u_0$ and $d_0$ satisfy (\ref{initden})-(\ref{bd}). Assume that $\rho_0\in C^1(\Omega)$, $u_0\in H^1(\Omega)$ and $d_0\in H^{2}(\Omega)$. Suppose that $(\rho, u, d)$ is a weak solution to the system \eqref{LCD} in Theorem \ref{exThm}.  Then\\
1. There is a positive small number $\e_0$ such that if
\begin{equation}\label{smalli}
\|u_0\|_{H^1(\Omega)}^2+\|\nabla d_0\|_{H^1(\Omega)}^2\leq\epsilon_0,
\end{equation}
then the system (\ref{LCD}) has a classical solution $(\rho,u,d)$ in the time period $(0, T)$, for all $T>0$. That is, for some $\alpha\in(0,1)$
\begin{equation}\label{reg}\begin{split}
u\in C^{1+\alpha/2,2+\alpha}((0,T)\times\Omega)\\
\nabla p\in C^{\alpha/2,\alpha}((0,T)\times\Omega)\\
d\in C^{1+\alpha/2,2+\alpha}((0,T)\times\Omega)\\
\rho\in C^1((0,T)\times\Omega).
\end{split}\end{equation}
2. For general data, there exists a positive number $\delta_0 = \delta_0(\rho_0, u_0, d_0)$ such that (\ref{reg}) holds in the interval $(0,T)$ for $T\leq\delta_0$.
\end{Theorem}

\begin{Corollary}\label{Cor2}Suppose in addition to the hypothesis in Theorem \ref{Mthm2}, that $\rho_0\in C^1(\bar\Omega)$,  $u_0\in C^{2+\alpha}(\bar\Omega)$ and $d_0\in C^{2+\alpha}(\bar\Omega)$, then the solution is regular up to the boundary for data small in the sense (\ref{smalli}) or for large data and sufficiently short time.
\end{Corollary}

The proof of the regularity of the solution to the system \eqref{LCD} in dimension 3 takes the same approach as in dimension 2 but is much more complicated. First in contrast to the cases of dimension 2, we only get the Ladyzhenskaya energy estimates when either the initial data is small in the sense as described in (\ref{smalli}) or $T$ is small. Our calculations and estimates are based on those in \cite{LL}, with interesting modifications. We use ideas of \cite{LL} making it work in a rather different way. We keep the potentially small terms $\|u\|_{L^2}$ and $\|\nabla d\|_{L^2}$ instead of throwing them away. This gives a more unified way to derive the Ladyzhenskaya energy estimates in the cases:
\begin{itemize}
\item of small data
\item for short time.
\end{itemize}
After having the Ladyzhenskaya energy estimates, in contrast to the two dimensional case, we do not have the H\"{o}lder continuity for the fluid density. Instead we observe that we have small oscillations of the density over small balls in $\Omega\times[0,T]$ provided that either the initial data is small or for short time. This turns out to be enough to carry out the frozen coefficient method to improve the regularity of the fluid velocity. We refer the reader to \cite{LS1} for a reference of the frozen coefficient method. We give the idea of the method in Appendix 6.
Our key lemma on the oscillation of the fluid density is as follows:

\begin{Lemma}\label{density} Suppose that $\Omega$ is a smooth bounded domain in $\R^3$. Let $\rho_0$, $u_0$ and $d_0$ satisfy (\ref{initden})-(\ref{bd}). Assume that $\rho_0\in C^1(\bar\Omega)$, $u_0\in H^1(\Omega)$ and $d_0\in H^{2}(\Omega)$. Suppose that $(\rho, u, d)$ is a weak solution to the system \eqref{LCD} in Theorem \ref{exThm}. Let $t_1\in(0,T)$ and $p\in\Omega$, define
$$
\Aa_{(p,t_1)}=(B_{r_0}(p)\cap\Omega)\times([t_1-r_0,t_1+r_0]\cap[0,T]).
$$
Then, for any $\epsilon>0$, there exists  $\epsilon_0>0$ and $r_0>0$ such that for $p\in \Omega$ and all $T>0$,
\bg\label{os}
\sup_{(q,t_2)\in\Aa_{(p,t_1)}}|\rho(q, t_2)- \rho(p, t_1)| \leq \epsilon,
\ed
provided that either
$$
\|u_0\|_{H^1(\Omega)}^2+\|\nabla d_0\|_{H^2(\Omega)}^2\leq\epsilon_0 \ \ \mbox { or } T\leq \delta_0.
$$
\end{Lemma}

\begin{Remark}
The interior regularity in Theorem \ref{Mthm1} and Theorem \ref{Mthm2} is obtained by bootstrapping argument.
\end{Remark}

To close this section, we state the uniqueness of solution in the following sense,
\begin{Theorem}\label{unique}
Let $(\rho,u,d)$ be the solution to system (\ref{LCD}) and (\ref{initden})-(\ref{bd}) obtained in Corollary \ref{Cor1} for two dimensions or in Corollary \ref{Cor2} for three dimensions.
Let $(\bar\rho,\bar u,\bar d)$ be a weak solution to system (\ref{LCD}) with (\ref{initden})-(\ref{bd}) satisfying the following energy inequalities:
\begin{align}\label{weakeng}
&\int_{\Omega}|\bar\rho|^2 dx \leq \int_{\Omega}|\rho_0|^2 dx,\\
&\int_{\Omega}\frac 1 2\bar\rho|\bar u|^2+\frac 1 2|\nabla\bar d|^2+
F(\bar d)dx+\int_0^T\int_{\Omega}|\nabla\bar u|^2+|\triangle\bar d-f(\bar d)|^2dxdt\\
&\leq \int_{\Omega}\frac 1 2\rho_0|u_0|^2+\frac 1 2|\nabla d_0|^2+F(d_0)dx.\notag
\end{align}
Then $(\rho,u,d)\equiv (\bar\rho,\bar u,\bar d)$.
\end{Theorem}
\begin{Remark} Note that the weak solutions established in \cite{Liu} and \cite{JT}  satisfy the inequalities recorded above. Hence these solutions coincide with solutions we obtained in Corollary \ref{Cor2} due to the uniqueness Theorem \ref{unique}.
\end{Remark}


\bigskip

\section{Classic solutions to Nematic Liquid Crystals Systems in $\R^2$}

In this section we use the Ladyzhenskaya energy method \cite{La} to show that $u \in L^\infty([0, T], H^1(\Omega))\bigcap L^2([0, T], H^2(\Omega))$ and $d\in L^\infty([0, T], H^2(\Omega))\bigcap L^2([0, T], H^3(\Omega))$. Then we apply a result from \cite{AKM} to get the H\"{o}lder continuity of the density $\rho$. With the continuity of the density $\rho$ we apply the so-called frozen coefficient method to get higher regularities for $\rho, p, u, d$ by bootstrapping. The proof of the Ladyzhenskaya energy inequality is similar to that in \cite{LL}. The key inequalities used often in this paper are the following Gagliardo-Nirenberg inequality (cf. \cite{Ev}) :

\begin{Lemma}\label{leGN} (Gagliardo-Nirenberg) If $\Omega$ is a smooth bounded domain in $\mathbb{R}^n$, then
\bg\label{GNd2+}
\|u\|_{L^4(\Omega)}^4\leq C\|u\|_{L^2(\Omega)}^2(\|\nabla u\|_{L^2(\Omega)}^2 + \|u\|_{L^2(\Omega)}^2),
\ed
when $n=2$ and
\bg\label{GNd3+}
\|u\|_{L^4(\Omega)}^4\leq C\|u\|_{L^2(\Omega)}(\|\nabla u\|_{L^2(\Omega)}^2 + \|u\|_{L^2(\Omega)}^2)^\frac 32,
\ed
when $n=3$. Moreover, if $u|_{\partial\Omega} = 0$, then
\bg\label{GNd2}
\|u\|_{L^4(\Omega)}^4\leq C\|u\|_{L^2(\Omega)}^2\|\nabla u\|_{L^2(\Omega)}^2,
\ed
when $n=2$ and
\bg\label{GNd3}
\|u\|_{L^4(\Omega)}^4\leq C\|u\|_{L^2(\Omega)}\|\nabla u\|_{L^2(\Omega)}^3,
\ed
when $n=3$.
\end{Lemma}

\subsection{The $L^\infty(H^1)$ and $L^2(H^2)$ estimates of the velocity}

Our strategy is the same as in \cite{LL}. We first establish the desired bounds for the Galerkin approximating  solutions $(\rho_m, u_m, d_m)$ in the sequence that one has from Galerkin method when proving the existence (cf. Section 5 of the appendix in this paper).  By passing to the weak limit we then obtain the desired bounds for the weak solution $(\rho, u, d)$. First we set
\bg\label{phi}
\Phi_m^2(t)=\|\nabla u^m\|_{L^2}^2+\|\triangle d^m\|_{L^2}^2.
\ed
Then to show that $u^m \in L^\infty( H^1)\bigcap L^2(H^2)$ and $d^m\in L^\infty( H^2)\bigcap L^2(H^3)$ we calculate $\frac d{dt} \Phi_m^2$.  Namely,

\begin{align}\notag
\frac{1}{2}\frac{d}{dt}\Phi_m^2(t)&=\int_{\Omega}\nabla u^m\cdot\nabla u_t^mdx+\int_{\Omega}\triangle d^m\cdot\triangle d_t^mdx\\
&= - \int_{\Omega}\rho^m|u^m_t|^2dx - \int_\Omega |\nabla\triangle d^m|^2dx\notag\\
&+\int_{\Omega}\nabla\triangle d^m\cdot\nabla(u^m\cdot\nabla d^m)-\triangle d^m\cdot\triangle(f(d^m))dx\notag\\
&+\int_{\Omega}-\rho^m(u^m\cdot\nabla u^m)u_t^m-u_t^m\nabla d^m\triangle d^mdx\notag\\
&= - \int_{\Omega}\rho^m|u^m_t|^2dx - \int_\Omega |\nabla\triangle d^m|^2dx\notag\\
& + I + II + II + IV,\notag
\end{align}
here we used the equation \eqref{aNSE}. We find
\begin{align}\notag
I & \leq \varepsilon\|\nabla\triangle d^m\|_{L^2}^2+C\|\nabla u^m\|_{L^4}^2\|\nabla d^m\|_{L^4}^2+C\|u^m\|_{L^4}^2\|\triangle d^m\|_{L^4}^2 + C \|u^m\|_{L^4}^2\notag\\
&\leq \varepsilon\|\triangle u^m\|_{L^2}^2 + 2\varepsilon\|\nabla\Delta d^m\|_{L^2}^2
+ C\Phi_m^4 + C\Phi_m^3+ C \Phi_m^2 + C\Phi_m,\notag
\end{align}
by Gagliardo-Nirenberg inequality, the basic energy estimate (\ref{energy}) and the fact that
$$
\|\nabla^2 d^m\|_{L^4} \leq C(\|\Delta d^m\|_{L^4} + \|d^m\|_{L^4}).
$$
And similarly,
\begin{align}\notag
II & \leq C\Phi_m^2 + C,\notag
\end{align}
\begin{align}\notag
III & \leq \varepsilon\|\triangle u^m\|_{L^2}^2+\varepsilon\|u_t^m\|_{L^2}^2 + C\Phi_m^4 + C\Phi_m^3,\notag
\end{align}
\begin{align}\notag
IV &\leq \varepsilon\|\nabla\triangle d^m\|_{L^2}^2+\varepsilon\|u_t^m\|_{L^2}^2+ C\Phi_m^4 + C \Phi_m^2.\notag
\end{align}
Moreover, from equation (\ref{aNSE}), we have
\bg\label{secu}
\|\triangle u^m\|^2_{L^2}\leq 2M_2 \|u_t^m\|_{L^2}^2+ 2\|\nabla\Delta d^m\|_{L^2}^2 + C\Phi_m^4 + C\Phi_m^3 + C\Phi_m^2
\ed
and
\bg\label{-secu}
M_1\|u^m_t\|^2_{L^2}\leq  2\|\Delta u^m\|_{L^2}^2+ 2\|\nabla\Delta d^m\|_{L^2}^2 + C\Phi_m^4 +C\Phi_m^3 +C\Phi_m^2.
\ed
Therefore, we arrive at
\bg\label{2energyin}
\frac{d}{dt}\Phi_m^2(t) + \|\Delta u^m\|_{L^2}^2 + \|\nabla\Delta d^m\|_{L^2}^2 \leq C\Phi_m^4(t) + C,
\ed
which implies, as in \cite{LL}, again in the light of the basic energy inequality \eqref{energy},
\bg\notag
\Phi_m^2(t) \leq (\Phi^2(0)+C) e^{C \int_0^T(\Phi_m^2(t)+1)dt} \leq (\Phi(0)^2+C) e^{CT +C}
\ed
where $\Phi^2(0)=\|\nabla u_0\|_{L^2}^2+\|\triangle d_0\|_{L^2}^2$, 
and
\bg\notag
\int_0^T\|\Delta u^m\|_{L^2}^2 dt + \int_0^T\|\nabla\Delta d^m\|_{L^2}^2dt \leq C.
\ed

Hence we have obtained

\begin{Theorem}\label{thm-phi}Suppose that $\Omega$ is a smooth bounded domain in $\mathbb{R}^2$ and that $\Phi(0)^2=\|\nabla u_0\|_{L^2}^2+\|\triangle d_0\|_{L^2}^2<\infty$. Suppose the initial data $(\rho_0,u_0,d_0)$ is as in Theorem \ref{Mthm1}.Then there exists a constant $C$ such that the approximating solutions $(\rho^m,u^m,d^m)$ of the system \eqref{LCD} obtained in Theorem \ref{exThm}  satisfies
\begin{equation}\label{phiest}
\|\nabla u^m\|_{L^2}^2 + \|\Delta d^m\|_{L^2}^2 \leq e^{CT +C} (\Phi^2(0)+C)
\end{equation}
for all $t\in [0, T]$ and
\bg\label{L^2H^2}
\int_0^T\|\Delta u^m\|_{L^2}^2 dt + \int_0^T\|\nabla\Delta d^m\|_{L^2}^2dt \leq C.
\ed
\end{Theorem}

From the last theorem it follows that
\begin{Corollary}Under the same hypothesis of last theorem, there exists a solution $(\rho, u, d)$ to (\ref{LCD}) which satisfies the energy inequalities
\begin{equation}\label{wphiest}
\|\nabla u\|_{L^2}^2 + \|\Delta d\|_{L^2}^2 \leq e^{CT +C} (\Phi^2(0)+C)
\end{equation}
for all $t\in [0, T]$ and
\bg\label{wL^2H^2}
\int_0^T\|\Delta u\|_{L^2}^2 dt + \int_0^T\|\nabla\Delta d\|_{L^2}^2dt \leq C.
\ed
 \end{Corollary}
\pf
It follows by extracting a subsequence of the Galerkin approximations $(\rho^m, u^m, d^m)$ and passing to the limit.

\subsection{H\"{o}lder continuity of the fluid density}

Next we  recall a regularity lemma for the transport equation from \cite{AKM} to get H\"{o}lder continuity for the fluid density $\rho$.

\begin{Lemma}\label{Le8}  (\cite{AKM})  Suppose that $u\in L^\infty([0, T], H^1(\Omega))\bigcap L^2([0, T], H^2(\Omega))$. And suppose that
$$
\rho_t + u\cdot\nabla\rho = 0
$$
in $\Omega\times (0, T)$. Assume that $\rho(0)\in C^1(\Omega)$ and that $\Omega\subset\R^2$ is smooth and bounded. Then $\rho\in C^\alpha(\Omega\times[0,T])$ for some $\alpha\in (0,1)$ which depends only on the initial data, $T$ and $\Omega$.
\end{Lemma}


\bigskip
\subsection{Proof of Theorem \ref{Mthm1}}

The proof is more or less standard, particularly after the work in \cite{AKM}. We sketch a proof here for completeness. Rewrite the third equation in (\ref{LCD}) as
\bg\label{red}
d_t-\triangle d=-u\cdot\nabla d-f(d).
\ed
It yields from the basic energy inequality (\ref{basiceng}) and estimate (\ref{phiest}) that ,
$$
u\cdot\nabla d\in L^\infty(0,T;L^r(\Omega)), \ \mbox { for any } r\in[1,\infty).
$$
Then due to standard estimates of solutions to parabolic equations (see \cite{LS} and \cite{Am}),
from the equation (\ref{red}) we have
$$
d\in W^{1, r}( W^{2,s}), \ \mbox { for any } r,s\in[1,\infty),
$$
which implies that
\bg\label{dholder}
d\in C^{\alpha/2,1+\alpha} ([0,T]\times\bar\Omega)
\ed
for any $\alpha\in (0, 1)$. Now we go to the second equation in (\ref{LCD}) and use the above estimates for $d$ to improve the estimate on $u$ via the frozen coefficient method \cite{LS1}, after we have the H\"{o}lder continuity for the fluid density $\rho$, as done in \cite{AKM} to derive that
$$
u \in W^{1,p}( W^{2, q}),
$$
for some $q > 1$ and any $p>1$. Therefore we have
$$
u\in C^{\alpha/2,\alpha}([0, T]\times\bar\Omega).
$$
Next going back to the third equation in \ref{LCD}, by the standard H\"{o}lder estimates for parabolic equations (see \cite{LS} Theorem 5.1 in Chapter VII) we have
$$
d \in C^{1+\frac {\alpha}{ 2},2+\alpha}((0, T)\times\Omega),
$$
for some $\alpha\in (0, 1)$. Therefore we are able to go back to the the second equation in (\ref{LCD}), again via frozen coefficient method \cite{LS1}, to derive that
$$
u \in C^{1+\frac {\alpha}{ 2},2+\alpha}((0, T)\times\Omega).
$$
Finally $\rho\in C^1((0, T)\times\Omega)$ follows from the regularity of $u$ and the regularity of the pressure $p$ follows easily from the regularity of $(\rho, u, d)$ similarly as in \cite{AKM}.
\cbdu
\begin{Remark}By bootstrapping argument, we can get higher regularity.
\end{Remark}
Proof of Corollary \ref{Cor1} follows by Krylov's Theorem 10.3.3 in \cite{Kr1}, which for convenience of reader, we recall here \cite{Kr1},
\begin{Theorem}For any $h\in C^{\alpha/2,\alpha}((0,T)\times\Omega)$ and boundary $g\in C^{1+\alpha/2,2+\alpha}((0,T)\times\Omega)$ there exists a unique function $u\in C^{1+\alpha/2,2+\alpha}([0,T]\times\bar\Omega)$ satisfying the equation $Lu-u_t=h$ in $(0,T)\times\Omega)$ and
equal $g$ on boundary of $(0,T)\times\Omega)$.
\end{Theorem}
We apply this to both the second and third equation in (\ref{LCD}). Applying to the second equation, $h=u\cdot\nabla d+f(d)$, $g=d_0$ and $L=\triangle$. Applying to the third equation, $h=u\cdot\nabla u+\nabla p+\nabla\cdot(\nabla d\otimes\nabla d)$, $g=0$ and $L=\triangle$. Proceeding as in the last theorem yields $h\in C^{\alpha/2,\alpha}((0,T)\times\Omega)$ for both equations. On the boundary of $\Omega$, $u=0$, $d=d_0$ and $d_0\in C^{2+\alpha}(\Omega)$ by hypothesis. Hence the conclusion of the Corollary follows.

\bigskip
\section{Classical Solution to Nematic Liquid Crystals System in 3D}

In this section we consider the solutions of (\ref{LCD}) for $\Omega\in\R^3$ a bounded domain. We establish regularity in two cases
\begin{itemize}
\item global regularity with small initial data
\item short time regularity.
\end{itemize}
First we adopt the idea from \cite{LL} to derive the Ladyzhenskaya energy estimates. In contrast to the cases of dimension 2, we will only get the Ladyzhenskaya energy estimates for the above two cases as in \cite{LL}. Our calculations and estimates are similar to those in \cite{LL}, with  interesting modifications. After having the Ladyzhenskaya energy estimates, unlike the case of dimension 2, we will not have H\"{o}lder continuity for the fluid density right away. Instead,
we observe that we can have the oscillation of the density over small balls in $\Omega\times[0,T]$ to be small, provided that either the initial data is small or we work for short time, which turns out to be enough to carry out the frozen coefficient method to improve the regularity of the fluid velocity.

\bigskip
\subsection{Ladyzhenskaya Energy Estimates}

 Our derivation of the Ladyzhenskaya energy estimates in \cite{LL} in dimension 3 is rather an interesting modification of the original Ladyzhenskaya's argument for the pure fluid systems. In \cite {LL} for the argument to work, it needs either the viscosity for the fluid to be very large or the time to be very short. We use the same idea, but, instead we assume the initial data to be small or the time to be short while the viscosity of the fluid is a fixed constant. For the convenience of the arguments in our context, without loss of generality we take the constant to be 1. Let us set as before:
\bg\label{phi}
\Phi_m^2(t)=\|\nabla u^m\|_{L^2}^2+\|\triangle d^m\|_{L^2}^2
\ed
Then, as in the previous section, we will first derive the Ladyzhenskaya energy estimates for the Galerkin approximate solutions $(\rho^m, u^m, d^m)$ and pass to the weak limit to obtain the Ladyzgenskaya energy estimates for the weak solutions $(\rho, u, d)$ to the system \eqref{LCD} as desired.

Again, Using $u^m|_{\partial\Omega}=d^m_t|_{\partial\Omega}=0$, and $\triangle d^m|_{\partial\Omega}=0$, by integration by parts, it still follows that,
\begin{align}\notag
\frac{1}{2}\frac{d}{dt}\Phi_m^2(t)&=\int_{\Omega}\nabla u^m\cdot\nabla u_t^mdx+\int_{\Omega}\triangle d^m\cdot\triangle d_t^mdx\\
&= - \int_{\Omega}\rho^m|u^m_t|^2dx - \int_\Omega |\nabla\triangle d^m|^2dx\notag\\
&+\int_{\Omega}\nabla\triangle d^m\cdot\nabla(u^m\cdot\nabla d^m)-\triangle d^m\cdot\triangle(f(d^m))dx\notag\\
&+\int_{\Omega}-\rho^m(u^m\cdot\nabla u^m)u_t^m-u_t^m\nabla d^m\triangle d^mdx\notag\\
&= - \int_{\Omega}\rho^m|u^m_t|^2dx - \int_\Omega |\nabla\triangle d^m|^2dx\notag\\
& + I + II + II + IV,\notag
\end{align}
We will proceed the same way as we did in dimension 2 except the Gagliardo-Nirenberg inequality is different in dimension 3 from that in dimension 2. More importantly we will keep the terms $\|u^m\|_{L^2}$ and $\|\nabla d^m\|_{L^2}$ whenever necessary. We may derive
$$
I \leq \varepsilon \|\Delta u^m\|_{L^2}^2 + \varepsilon\|\nabla\Delta d^m\|_{L^2}^2 +   (\|u^m\|_{L^2} + \|\nabla d^m\|_{L^2}) (C\Phi_m^8 +C),
$$
$$
II \leq C\Phi_m^2 + C\|\nabla d^m\|_{L^2},
$$
$$
III \leq \varepsilon \|u_t^m\|_{L^2}^2 + \varepsilon\|\Delta u^m\|_{L^2}^2 + C\|u^m\|_{L^2}^2\Phi_m^8 + C\Phi_m^2,
$$
and
$$
IV\leq \varepsilon \|u_t^m\|_{L^2}^2 + \varepsilon\|\nabla\Delta d^m\|_{L^2}^2
+ C\|\nabla d^m\|_{L^2}^2\Phi_m^8 + C\|\nabla d^m\|_{L^2}^2.
$$
To relate $\Delta u^m$ back to $u_t^m$, from equation (\ref{aNSE}), we have that
\begin{align}
\|\triangle u^m\|_{L^2}^2 &\leq 2M_2\|u_t^m\|_{L^2}^2 + 2\|\nabla\Delta d^m\|_{L^2}^2 + (\|u^m\|_{L^2}^2 + \|\nabla d^m\|_{L^2}^2)(C\Phi_m^8 +C) + C\Phi_m^2 \notag
\end{align}
and
\begin{align}
M_1\| u_t^m\|_{L^2}^2 &\leq 2\|\Delta u^m\|_{L^2}^2 + 2\|\nabla\Delta d^m\|_{L^2}^2 + (\|u^m\|_{L^2}^2 + \|\nabla d^m\|_{L^2}^2)(C\Phi_m^8 +C) + C\Phi_m^2. \notag
\end{align}
Therefore we arrive at
$$
\|\Delta u^m\|_{L^2}^2 + \|\nabla\Delta d^m\|_{L^2}^2 + \frac d{dt}\Phi_m^2 \leq C\Phi_m^2
+ (\|u^m\|_{L^2} + \|\nabla d^m\|_{L^2})(C\Phi_m^8 +C).
$$
Set
$$
\tilde\Phi_m^2 = \Phi_m^2+\|u^m\|_{L^2} + \|\nabla d^m\|_{L^2}
$$
and observe that
\begin{equation}\label{sm-ini}
\frac d{dt} \tilde\Phi_m^2 \leq (C + C(\|u^m\|_{L^2} + \|\nabla d^m\|_{L^2})\tilde\Phi_m^6)\tilde\Phi_m^2
\end{equation}
or set
$$
\tilde\Phi_m^2 = \Phi_m^2+\|u^m\|_{L^2} + \|\nabla d^m\|_{L^2}+1
$$
we have
\begin{equation}\label{sm-time}
\frac d{dt} \tilde\Phi_m^2 \leq C \tilde\Phi_m^8.
\end{equation}
It is clear from \eqref{sm-time} one can prove the Ladyzhenskaya energy estimates when $T$ is small.  To get the Ladyzhenskaya energy estimates for small initial data we will use an idea similar to that in \cite{LL}. Suppose that
$$
\|u_0\|_{H^1} ^2+ \|d_0\|_{H^2}^2 = \theta_0.
$$
Recall that by the basic energy estimate we have
$$
\|u^m\|_{L^2} + \|\nabla d^m\|_{L^2}\leq\|u_0\|_{H^1} ^2+ \|d_0\|_{H^2}^2 \leq\theta_0.
$$
Then we claim that, if $\theta_0$ is so small that
\begin{equation}\label{small}
C\theta_0 (4e^{C+1}\theta_0)^3 < 1,
\end{equation}
then
$$
\tilde\Phi_m^2 \leq 4 e^{C+1}\theta_0,
$$
for all $t$. First we prove the claim for $t\in [0, 1]$. Assume otherwise, there must be $t_0\in (0, 1)$ such that
\begin{equation}\label{wrong}
\tilde\Phi_m^2(t_0) = 4 e^{C+1}\theta_0
\end{equation}
and
$$
\tilde\Phi_m^2(t) \leq 4 e^{C+1}\theta_0
$$
for all $t\in (0, t_0]$. Therefore, from (\ref{sm-ini}), by the choice of $\theta_0$ in \eqref{small}, we have
$$
\frac d{dt} \tilde\Phi_m^2 \leq (C+1)\tilde\Phi_m^2
$$
for all $t\in (0, t_0)$ and $\tilde\Phi_m^2(0) \leq 2\theta_0$, which implies that
$$
\tilde\Phi_m^2(t_0) \leq e^{C+1}\tilde\Phi_m^2(0)\leq 2e^{C+1}\theta_0
$$
and thus contradicts \eqref{wrong}. For $t >1$, we simply observe, as in \cite{LL}, that
the basic energy inequality \eqref{basiceng}
$$
\int_{t-1}^t \tilde\Phi^2 _m(t)dt \leq 2\theta_0
$$
implies that there is $t_0\in (t-1, t)$ such that
$$
\tilde\Phi_m^2(t_0) \leq 2\theta_0.
$$
Then one may repeat the above argument to conclude that
$$
\tilde\Phi_m^2(t)\leq 4 e^{C+1}\theta_0.
$$
As in the case of 2 dimension, passing to the limit the existence of weak solutions will follow from the uniform estimates we have obtained. This weak solution satisfies

\begin{Theorem}\label{lady-en} Suppose that $\Omega$ is a smooth bounded domain in $\R^3$. Let $\rho_0$, $u_0$ and $d_0$ satisfy (\ref{initden})-(\ref{bd}). Assume that
$$
\|u_0\|_{H^1} + \|d_0\|_{H^2} < \infty.
$$
Let $(\rho, u, d)$ be a weak solution to the system \eqref{LCD} with data $(\rho_0,u_0,d_0)$. There is $\epsilon_0>0$ such that if
$$
\|u_0\|_{H^1} + \|d_0\|_{H^2} \leq \epsilon_0,
$$
then
\begin{equation}\label{lady3}
\|\nabla u\|_{L^2}^2 + \|\Delta d\|_{L^2}^2 + \int_0^T(\|\Delta u\|_{L^2}^2 + \|\nabla\Delta d\|_{L^2}^2)dt \leq C(\|u_0\|_{H^1} + \|d_0\|_{H^2})
\end{equation}
holds for any positive  $T$.
For data with no smallness condition there is a $\delta_0$ depending on the data such that (\ref{lady3}) holds for all
$T\leq \delta_0$.
\end{Theorem}

\bigskip

\subsection{Oscillation of the Fluid Density}

In this subsection, we show that the oscillation of the density $\rho$ over small balls in space time. This estimate follows by the estimate on $u$ from the Ladyzhenskaya energy estimates obtained in the previous subsection. For convenience of the readers we restate Lemma \ref{density}.

\begin{Lemma} Suppose that $\Omega$ is a smooth bounded domain in $\R^3$. Let $\rho_0$, $u_0$ and $d_0$ satisfy (\ref{initden})-(\ref{bd}). Assume that $\rho_0\in C^1(\bar\Omega)$, $u_0\in H^1(\Omega)$ and $d_0\in H^{2}(\Omega)$. Suppose that $(\rho, u, d)$ is a weak solution to the system \eqref{LCD} constructed in Theorem \ref{exThm} ( see Appendix). Let $t_1\in(0,T)$ and $p\in\Omega$, define
$$
\Aa_{(p,t_1)}=(B_{r_0}(p)\cap\Omega)\times([t_1-r_0,t_1+r_0]\cap[0,T]).
$$
Then, for any $\epsilon>0$, there exists  $\epsilon_0>0$ and $r_0>0$ such that for $p\in \Omega$ and all $T>0$,
\bg\notag
\sup_{(q,t_2)\in\Aa_{(p,t_1)}}|\rho(q, t_2)- \rho(p, t_1)| \leq \epsilon,
\ed
provided that either
$$
\|u_0\|_{H^1(\Omega)}^2+\|\nabla d_0\|_{H^2(\Omega)}^2\leq\epsilon_0 \ \ \mbox { or } \ \  T\leq \delta_0.
$$
\end{Lemma}

\pf The equation of the conservation of mass is
\begin{equation}\label{mass}
\rho_t  + u\cdot \nabla \rho = 0.
\end{equation}
For reasonably regular velocity we may solve \eqref{mass} by the method of characteristics. Due to Osgood theorem \cite{Os}, since $u\in H^2$ for almost all $t\in [0, T]$ and in light of the Ladyzhenskaya energy estimates in Theorem \ref{lady-en}, there is a unique solution to the Cauchy problem corresponding to (\ref{mass}) by finding  trajectories of the liquid particles
\bg\notag
\frac{dy}{d\tau}=u(y,\tau), \ \ y|_{\tau=t}=x; \ \ x\in\Omega,\ \ t\in(0,T)
\ed
and defining
\begin{equation}\label{rho-sol}
\rho(x, t) = \rho_0(y(0, x, t)).
\end{equation}

\noindent
{\bf{Step 1}} \quad Fix a time $t\in [0, T]$. Let $x_1$ and $x_2$ be two arbitrary points from $\Omega$ satisfying $|x_1-x_2|\leq d<1$. For any $\tau\in(0,t)$, assume
$$
y_1=y(\tau,x_1,t), \ \ y_2=y(\tau,x_2,t),
$$
then the difference $z(\tau)=y_1-y_2$ is the solution to the Cauchy problem
$$
\frac{dz}{d\tau}=u(y_1,\tau)-u(y_2,\tau), \ \ z|_{\tau=t}=x_1-x_2.
$$
For $0\leq\alpha\leq 1/2$, by a standard Sobolev embedding theorem in \cite{Ev}, we have
$$
\frac{d|z|}{d\tau}\leq C\|u(\tau)\|_{H^2(\Omega)}|z|^\alpha, \ \ \tau\in(0,t).
$$
Integrating form $\tau$ to $t$,
\bg\label{z}
|z|^{1-\alpha}\leq |x_1-x_2|^{1-\alpha}+CT^{\frac{1}{2}}\|u\|_{L^2(0,T;H^2(\Omega))}.
\ed
As a consequence of Theorem \ref{lady-en}, we know that
$|z|$ is as small as needed provided that $|x_1-x_2|$ and $\|u_0\|_{H^1} + \|d_0\|_{H^2}$ are small or $T$ is small. Therefore
$$
|\rho (x_1, t) - \rho(x_2, t)| = |\rho_0(y(0, x_1, t))  - \rho_0(y(0, x_2, t))| \leq C|z|
$$
is as small as one wants provided that $|x_1-x_2|$ and $\|u_0\|_{H^1} + \|d_0\|_{H^2}$ are small or $T$ is small.

\vskip 0.1in\noindent
{\bf{Step 2}} \quad Fix a point $x\in \Omega$. Let $t_1, t_2\in[0,T]$ arbitrary, let
$$
y_1=y(\tau,x,t_1), \ \ y_2=y(\tau,x,t_2).
$$
Assume that $x'=y_2|_{\tau=t_1}$. Then due to uniqueness, the integral curve $y_2(\tau)$ can be considered as a solution to the Cauchy problem
$$
\frac{dy_2}{d\tau}=u(y_2,\tau), \ \ y_2|_{\tau=t_1}=x'
$$
with initial data at $\tau=t_1$. Hence, the difference $z(\tau)=y_1-y_2$ is the solution of the Cauchy problem
$$
\frac{dz}{d\tau}=u(y_1,\tau)-u(y_2,\tau), \ \ z|_{\tau=t_1}=x-x'.
$$
By the definition of $x'$, we have
$$
x'=x-\int_{t_1}^{t_2}u(y_2(s,x,t_2),s)ds.
$$
Therefore, again due to a standard Sobolev embedding Theorem,
\begin{align}\notag
|x-x'| \leq \int_{t_1}^{t_2}|u(y,s)|ds
\leq C \int_{t_1}^{t_2} \|u\|_{H^2(\Omega)} ds\leq  C \|u\|_{L^2(0,T;H^2(\Omega))} |t_1-t_2|^\frac 12.\notag
\end{align}
By Step 1, we conclude that
$$
|\rho(x, t_1) - \rho(x, t_2)| = |\rho_0(y(0, x, t_1))- \rho_0(y(0, x, t_2))| = |\rho_0(y(0, x, t_1)) - \rho_0(y(0, x', t_1))|
$$
is as small as needed provided that $|t_1-t_2|$ and $\|u_0\|_{H^1} + \|d_0\|_{H^2}$ are small or $T$ is small.
This completes the proof of the Lemma.
\cbdu

\bigskip

\subsection{Proof of Theorem \ref{Mthm2}}

In this subsection, we finish the proof of Theorem \ref{Mthm2}. This is done by combining Theorem \ref{lady-en}, Lemma \ref{density} from the previous subsections, the frozen coefficient technique applied to $L^p(L^q)$ estimates, H\"{o}lder estimates and a bootstrapping argument between the three equations of the system \eqref{LCD}. \\
In the appendix we briefly describe the frozen coefficient method. For simplicity, we show how  the method works for the density dependent NSE and obtain the estimate for $L^q$ space with $q>\sqrt 3$. We apply now this method  to our approximating solutions $(\rho^m, u^m, d^m)$. These solutions for each $m$ satisfy the conclusion of Theorem \ref{Mthm2}. We now show that the estimates are uniformly in $m$ applying the frozen coefficient method, in which case we let $q=2$. That allows us to pass to the limit $(\rho, u, d)$ and obtain  (\ref{reg}) for weak solutions. We now show briefly the steps to yield the necessary uniform bound.\\
Here to simplify the notation we denote the approximating solutions by $(\rho, u, d)$.
We notice first that since
$$
u\in L^\infty(0,T;L^6(\Omega)), \ \ \nabla d\in L^\infty(0,T;L^6(\Omega)),
$$
we have
$$
u\cdot\nabla d\in L^\infty(0,T;L^3(\Omega)).
$$
By standard parabolic estimates on the third equation in (\ref{LCD}) (cf. \cite{LS} and \cite{Am}), we have
$$
d\in W^{1,p}(W^{2,3}),
$$
for all $1 < p <\infty$, which implies that $\nabla d\in L^\infty(0,T;L^q(\Omega))$ for any $q\in(1,\infty)$. Thus, we have
$$
u\cdot\nabla d\in L^\infty(0,T;L^q(\Omega)), \ \ \forall q\in(1,6).
$$
Applying the same standard parabolic estimate on the third equation in (\ref{LCD}) yields
\bg\label{dw2r}
d\in W^{1, p}(W^{2,q}), \forall p\in (1, \infty) \ \text{and} \ q \in(1,6),
\ed
which implies that $d\in C^{\alpha/2,1+\alpha}([0,T]\times\bar\Omega)$ for some $\alpha \in (0, 1)$ and
$$
\nabla d\Delta d \in L^\infty(0, T; L^q(\Omega)), \quad \forall q \in (1, 6).
$$

In the second equation of (\ref{LCD}), the estimates for the conservation of momentum with constant density can be extended to the non-constant density cases when Lemma \ref{density} is available. This is done via the frozen coefficient method.

We know that
$$
u\cdot\nabla u\in L^\infty(0,T;L^{\frac{3}{2}}(\Omega)).
$$
Now we apply the frozen coefficient method using the oscillation estimates for the density, (Lemma \ref{density})  to yield
$$
u\in W^{1, p}(W^{2,3/2}), \quad \forall p\in (1, \infty)
$$
Thus $u\in L^\infty(0,T;W^{1,3}(\Omega))$ and $u\cdot\nabla u\in L^\infty(0,T;L^2(\Omega))$. Repeating the above argument yields
\bg\label{uh2}
u\in W^{1, p}(W^{2,2}), \quad \forall p\in (1, \infty),
\ed
from where it follows that $u\in C^\alpha(\Omega\times[0,T])$.

Back to the third equation in (\ref{LCD}),  we conclude that
$$
d \in C^{1+\frac \alpha 2,2+\alpha}((0, T)\times\Omega)
$$
for some $\alpha\in (0, 1)$ via the standard H\"{o}lder estimates for parabolic equations (cf. \cite{LS}, Chapter VII).  From here the argument for the pure fluid systems works with no further significant modifications. This completes the proof of Theorem \ref{Mthm2}.\\

Proof of Corollary \ref{Cor2}: Follows by Krylov's Theorem \cite{Kr2} just as in the case of two dimension.

\bigskip

\section{Uniqueness of Solution}
In this section we establish Theorem \ref{unique}. For  the LCD system with constant density, Lin and Liu \cite{LL} proved that the solution $(u,d)$ is unique provided $u\in L^\infty(0,T;H^1)$ and $d\in L^\infty(0,T;H^2)$. The idea is to calculate the energy law satisfied by the difference of two solutions and establish a Gronwall's inequality. In our case, to calculate the energy law of the difference of two solutions it has some extra terms involving the density. Hence the estimates are more involved  requiring  additional  bounds on the strong solution $(\rho,u,d)$ to yield a Gronwall's inequality. In 2D, we need
\bg\label{more2}
\nabla\rho,\nabla u\in L^\infty((0,T)\times\Omega), \ \ u_t, u\cdot\nabla u\in L^\infty(0,T;L^q(\Omega)), q>2.
\ed
In 3D, we need
\bg\label{more3}
\nabla\rho,\nabla u\in L^\infty((0,T)\times\Omega), \ \ u_t, u\cdot\nabla u\in L^\infty(0,T;L^3(\Omega)).
\ed
With the assumption on data, $\rho_0\in C^1(\bar\Omega)$,  $u_0\in C^{2+\alpha}(\bar\Omega)$ and $d_0\in C^{2+\alpha}(\bar\Omega)$,  the solution $(\rho,u,d)$ from Corollary \ref{Cor1} or Corollary  \ref{Cor2} satisfies (\ref{more2}) or (\ref{more3}), respectively .\\


Proof of Theorem \ref{unique}: First recall that the solution $(\rho, u,d)$ from Theorem \ref{Mthm1} and Theorem \ref{Mthm2} satisfies energy equality:
\begin{align}\label{engeq}
&\int_{\Omega}\frac 1 2\rho|u|^2+\frac 1 2|\nabla d|^2+
F(d)dx+\int_0^T\int_{\Omega}|\nabla u|^2+|\triangle d-f(d)|^2dxdt\\
&= \int_{\Omega}\frac 1 2\rho_0|u_0|^2+\frac 1 2|\nabla d_0|^2+F(d_0)dx.\notag
\end{align}
The density   $\rho$ is the strong solution of  the transport equation, hence it  satisfies that
$$
\int_\Omega\rho^2dx=\int_\Omega\rho_0^2dx.
$$
On the other hand side,  $\bar\rho$ is a weak solution of the transport equation and $M_1\leq\bar\rho\leq M_2$.   We have by hypothesis that
\bg\label{barrho}
\int_\Omega\bar\rho^2dx\leq\int_\Omega\rho_0^2dx.
\ed
Thus,
\begin{align}\label{rhodiff}
\frac 1 2\int_\Omega|\rho-\bar\rho|^2dx&=\frac 1 2\int_\Omega\rho^2dx+\frac 1 2\int_\Omega\bar\rho^2dx-\int_\Omega\rho\bar\rho dx\\
&\leq\int_\Omega\rho_0^2dx-\int_\Omega\rho\bar\rho dx.\notag
\end{align}
Since $\rho\in C^1([0,T]\times\bar\Omega)$,  we can take $\rho$ as a test function. Thus, multiplying
\bg\notag
\bar\rho_t+\bar u\cdot\nabla\bar\rho=0
\ed
by $\rho$
and integrating by parts yields
\begin{align}\label{limit}
\int_\Omega\rho_0^2dx-\int_\Omega\rho\bar\rho dx&=-\int_0^t\int_\Omega\bar\rho\rho_t dsdx-\int_0^t\int_\Omega(\bar u\cdot\nabla\rho)\bar\rho dsdx\\
&=\int_0^t\int_\Omega(u\cdot\nabla\rho)\bar\rho dsdx-\int_0^t\int_\Omega(\bar u\cdot\nabla\rho)\bar\rho dsdx.\notag
\end{align}
Here we used again that $\rho$ is a classical solution of the transport equation. Substituting the equality (\ref{limit}) in (\ref{rhodiff}) gives
\begin{align}\label{rhodiff1}
\frac 1 2\int_\Omega|\rho-\bar\rho|^2dx&\leq\int_0^t\int_\Omega\bar\rho(u-\bar u)\nabla\rho dsdx\\
&=\int_0^t\int_\Omega(\rho-\bar\rho)(u-\bar u)\nabla\rho dsdx.\notag
\end{align}
Next, calculate the following term
\begin{align}\notag
&\frac 1 2\int_\Omega\bar\rho|u-\bar u|^2dx+\frac 1 2\int_\Omega|\nabla d-\nabla\bar d|^2dx\\
&=\frac 1 2\int_\Omega(\bar\rho-\rho)|u|^2dx+\frac 1 2\int_\Omega\rho|u|^2dx+\frac 1 2\int_\Omega\bar\rho|\bar u|^2dx-\int_\Omega\bar\rho u\cdot\bar udx\notag\\
&+\frac 1 2\int_\Omega|\nabla d|^2dx+\frac 1 2\int_\Omega|\nabla\bar d|^2dx-\int_\Omega \nabla d\otimes\nabla\bar ddx\notag,
\end{align}
where $\nabla d\otimes\nabla \bar d$ denotes the $3\times 3$ matrix whose $ij$-th entry is given by $\nabla_i d\cdot\nabla_j\bar d$ for $1\leq i,j\leq 3$.\\
Using energy equality (\ref{engeq}) for the regular  solution $(\rho,u,d)$ and inequality (\ref{weakeng})  for the weak  solution $(\bar{\rho},\bar{ u},\bar d)$ combined with the last equality yields
\begin{align}\label{uddiff}
&\frac 1 2\int_\Omega\bar\rho|u-\bar u|^2dx+\frac 1 2\int_\Omega|\nabla d-\nabla\bar d|^2dx\\
&\leq-\int_0^t\int_\Omega|\nabla u-\nabla\bar u|^2dxdt-\int_0^t\int_\Omega|\triangle d-\triangle\bar d|^2dxdt\notag\\
&+\frac 1 2\int_\Omega(\bar\rho-\rho)|u|^2dx-\int_\Omega\bar\rho u\otimes\bar u-\rho_0|u_0|^2dx-\int_\Omega \nabla d\otimes\nabla\bar d-|\nabla d_0|^2 dx\notag\\
&-\int_\Omega F(d)dx-\int_\Omega F(\bar d)dx+2\int_\Omega F(d_0)dx-\int_0^t\int_\Omega |f(d)|^2dxdt\notag\\
&-\int_0^t\int_\Omega |f(\bar d)|^2dxdt-2\int_0^t\int_\Omega\nabla u\cdot\nabla\bar udxdt-2\int_0^t\int_\Omega\triangle d\cdot\triangle\bar ddxdt\notag\\
&+2\int_0^t\int_\Omega\triangle d\cdot f(d)dxdt+2\int_0^t\int_\Omega\triangle\bar d\cdot f(\bar d)dxdt.\notag
\end{align}
Since $u,d\in C^{1+\alpha/2,2+\alpha}([0,T]\times\bar\Omega)$, we can take $u,d$ as test functions for the weak
solution $\bar u,\bar d$. Thus it follows that
\begin{align}\label{uweak}
\int_\Omega\bar\rho u\otimes\bar u-\rho_0|u_0|^2dx&=\int_0^t\int_\Omega\bar\rho\bar uu_tdxdt+\int_0^t\int_\Omega\bar\rho\bar u(\bar u\cdot\nabla u)dxdt\\
&-\int_0^t\int_\Omega\bar\rho\bar u(\bar u\cdot\nabla u)dxdt-\int_0^t\int_\Omega\bar\rho\bar u(\bar u\cdot\nabla u)dxdt,\notag
\end{align}
\begin{align}\label{dweak}
\int_\Omega \nabla d\otimes\nabla\bar d-|\nabla d_0|^2 dx&=-\int_0^T\int_\Omega\triangle\bar dd_tdxdt+\int_0^T\int_\Omega\triangle d(\bar u\cdot\nabla\bar d)dxdt\\
&-\int_0^T\int_\Omega\triangle\bar d\cdot\triangle ddxdt+\int_0^T\int_\Omega f(\bar d)\triangle ddxdt.\notag
\end{align}
Indeed, formally, to get (\ref{uweak}), we multiply equation
$$
(\bar\rho\bar u)_t+\div(\bar\rho\bar u\otimes\bar u)+\nabla\bar p=\triangle\bar u-\nabla\cdot(\nabla\bar d\otimes\nabla\bar d)
$$
by $u$
and integrate by parts. To get (\ref{dweak}), we multiply equation
$$
\bar d_t+\bar u\cdot\nabla\bar d=\triangle\bar d-f(\bar d)
$$
by $\Delta d$ and integrate by parts.\\

Substitute (\ref{uweak}) and (\ref{dweak}) in (\ref{uddiff}). Add (\ref{rhodiff1}) combined with the equalities in equations  (\ref{LCD}), yields
\begin{align}\label{engddiff}
&\frac 1 2\int_\Omega|\rho-\bar\rho|^2dx+\frac 1 2\int_\Omega\bar\rho|u-\bar u|^2dx+\frac 1 2\int_\Omega|\nabla d-\nabla\bar d|^2dx\\
&\leq-\int_0^t\int_\Omega|\nabla u-\nabla\bar u|^2dxdt-\int_0^t\int_\Omega|\triangle d-\triangle\bar d|^2dxdt\notag\\
&+\int_0^t\int_\Omega(\bar\rho-\rho)(u-\bar u)\nabla\rho dxdt-\int_0^t\int_\Omega\bar\rho\nabla u|u-\bar u|^2dx\notag\\
&+\int_0^t\int_\Omega (\bar\rho-\rho)(u-\bar u)(u_t+u\cdot\nabla u)dxdt+\int_0^t\int_\Omega u(\nabla d-\nabla\bar d)(\triangle d-\triangle\bar d)dxdt\notag\\
&-\int_0^t\int_\Omega (u-\bar u)(\nabla d-\nabla\bar d)\triangle ddxdt
+\int_0^t\int_\Omega(\triangle d-\triangle\bar d)(f(d)-f(\bar d))dxdt.\notag
\end{align}
Recall that the regular solution $(\rho,u,d)$ satisfies (\ref{more2}) in 2D and (\ref{more3}) in 3D.  By H\"older  and Gagliardo-Nirenberg inequalities on the right hand side of (\ref{engddiff}), it follows that
\begin{align}\label{Grnw}
&\frac 1 2\int_\Omega|\rho(t)-\bar\rho(t)|^2+\bar\rho(t)|u(t)-\bar u(t)|^2+|\nabla d(t)-\nabla\bar d(t)|^2dx\\
&\leq C\int_0^t\int_\Omega|\rho-\bar\rho|^2+|u-\bar u|^2+|\nabla d-\nabla\bar d|^2dxdt.\notag
\end{align}
To handle the last integral in (\ref{engddiff}), we used the fact that $|d|,|\bar d|\leq 1$ and hence $|f(d)-f(\bar d)|\leq C|d-\bar d|$ by the definition of $f(d)$. Thus,
\bg\notag
\int_0^t\int_\Omega|f(d)-f(\bar d)|^2dxdt\leq C\int_0^t\int_\Omega|d-\bar d|^2dxdt\leq C(\Omega)\int_0^t\int_\Omega|\nabla d-\nabla\bar d|^2dxdt
\ed
where the constant $C$ depends on space domain $\Omega$ not on time $T$, and $C$ depends on the dimension of the space.\\

Using the lower bound of $\bar\rho\geq M_1>0$, and Gronwall's inequality to (\ref{Grnw}) yields
\begin{align}\notag
&\frac 1 2\int_\Omega|\rho(t)-\bar\rho(t)|^2+\bar\rho|u(t)-\bar u(t)|^2+|\nabla d(t)-\nabla\bar d(t)|^2dx\\
&\leq \int_\Omega|\rho(0)-\bar\rho(0)|^2+\bar\rho(0)|u(0)-\bar u(0)|^2+|\nabla d(0)-\nabla\bar d(0)|^2dxe^{Ct}\notag\\
&=0\notag
\end{align}
for all $t>0$ which implies
$$
\bar\rho-\rho=\bar u-u=\bar d-d\equiv 0.
$$
This completes the proof of Theorem \ref{unique}.

\bigskip

\section{Appendix: Existence of Weak Solutions}

In this appendix we sketch an existence theorem for Galerkin approximations.
As mentioned in the introduction the existence of global weak solutions to the flow of nematic liquid crystals have been established in \cite{Liu} and in \cite{JT} for non constant fluid density. Unfortunately the Ladyzhenskaya energy estimates do not seem to work for the Galerkin approximate solutions constructed in \cite{Liu} and in \cite{JT}.  The Galerkin approximate solutions constructed here will possess the Ladyzhenskaya energy estimates. When the initial fluid density has a positive lower bound, we are able to derive estimates on $u^m_t$ and $d^m_t$ so that we can employ the compactness lemma of Lions-Aubin. Since the Galerkin method has been widely used for fluid systems as well as on the system of liquid crystals we will be rather brief (cf. \cite{AKM},  \cite{LL},  \cite{Liu}, \cite{JT}).

\bigskip
\subsection{Galerkin Approximate Solutions}
We construct a sequence of Galerkin approximating solutions that are uniformly bounded. These are the solutions used in Section 2 and Section 3. The bounds obtained there through Ladyzhenskaya method yield a subsequence that will converge to the classical solution.\\
Let us first state the existence theorem for global weak solutions:

\begin{Theorem}\label{exThm}
Assume that $u_0\in L^2(\Omega)$ and $d_0\in H^1(\Omega)$ with $d_0|_{\partial\Omega}\in H^{3/2}(\Omega)$.
The system (\ref{LCD}) with initial boundary conditions (\ref{initden})-(\ref{bd}) has a weak solution $(\rho^m, u^m, d^m)$ for each $m=1,2,3, ...$ satisfying, for all $T\in(0,\infty)$,
\bg
0< M_1 \leq\rho^m\leq M_2
\ed
\bg
u^m\in L^2(0,T; H_0^1(\Omega))\cap L^\infty(0,T;L^2(\Omega))
\ed
\bg
d^m\in L^2(0,T; H^2(\Omega))\cap L^\infty(0,T;H^1(\Omega)), \quad |d^m| \leq 1,
\ed
and the energy inequality
\begin{equation}\label{basiceng}
\int_{\Omega}(|\triangle d^m-f(d^m)|^2+|\nabla u^m|^2)dx + \frac{d}{dt}\int_{\Omega}(\frac{1}{2}\rho |u^m|^2
+\frac{1}{2}|\nabla d^m|^2+F(d^m))dx\leq 0.
\end{equation}
\end{Theorem}

\proof
Let
$$
\Hh(\Omega)= \mbox { closure of } \{f\in C^\infty_0(\Omega,\mathbb{R}^3): \nabla\cdot f = 0\}   \mbox { in } L^2(\Omega,\mathbb{R}^3)
$$
and
$\left\{\phi_i\right\}_{i=1}^\infty$ be an orthonormal basis of $\Hh$ and satisfying:
\bg\notag
-\triangle\phi_i+\nabla P_i=\lambda_i\phi_i \mbox { in } \Omega
\ed
\bg\notag
\phi_i=0 \mbox { on } \partial\Omega
\ed
for $i=1,2,...$. Here $0<\lambda_1\leq\lambda_2...\leq\lambda_n\leq...$ with $\lambda_i\to\infty$. In other words, we choose an orthonormal basis of $\Hh$ which consists of the eigenfunctions of Stokes operator on $\Omega$ with vanishing Dirichlet boundary condition (see \cite{Tem}). Let
$$
P_m: \Hh\to \Hh_m=span\left\{ \phi_1,...\phi_m\right\}
$$
be the orthonormal projection. We seek approximate solutions $(\rho^m, u^m, d^m)$ with $u^m\in \Hh_m$, satisfy the following equations:
\bg\label{aden}
\rho_t^m+u^m\cdot\nabla\rho^m=0,\ \ \rho^m(x,0)=\rho_0(x)
\ed
\bg\label{aNSE}
P_m(\rho^m\frac{\partial}{\partial t}u^m)=P_m(\triangle u^m-\rho^mu^m\cdot\nabla u^m-\nabla\cdot(\nabla d^m\otimes\nabla d^m))
\ed
\bg\label{ad}
d_t^m+u^m\cdot\nabla d^m=\triangle d^m-f(d^m)
\ed
\bg\label{auinit}
u^m(x,0)=P_m u_0(x), \quad
d^m(x,0)=d_0(x), \quad d^m(x,t)|_{\partial\Omega}=d_0(x)|_{\partial\Omega}.
\ed
Let
\bg\label{au}
u^m(x,t)=\sum_{l=1}^mg_i^m(t)\phi_i(x),
\ed
with  $g_i^m(t)\in C^1[0,T]$. Hence (\ref{aNSE}) is equivalent to the following system of ordinary differential equations:
\bg\label{ODE}
\sum_{i=1}^m\Aa_i^{mj}(t)\frac{d}{dt}g_i^m(t)=-\sum_{i,k}^m\Bb_{ik}^{mj}(t)g_i^m(t)g_k^m(t)
-\sum_{i=1}^m\Cc_i^jg_i^m(t)+\Dd^{mj}(t),
\ed
for $j=1,2,..., m$, where
\begin{equation}\label{cff}
\begin{cases}
\Aa_i^{mj}(t)=\int_\Omega\rho^m(t)\phi_i(x)\phi_j(x)dx,\\
\Bb_{ik}^{mj}(t)=\int_\Omega\rho^m(t)(\phi_i(x)\cdot\nabla\phi_k(x))\phi_j(x)dx,\\
\Cc_i^{j}=\int_\Omega\nabla\phi_i(x)\cdot\nabla\phi_j(x)dx,\\
\Dd^{mj}(t)=\int_\Omega\sum_{k,l}(\frac{\partial}{\partial x_k}d^m\cdot\frac{\partial}{\partial x_l}d^m)\frac{\partial}{\partial x_l}\phi_j^k(x)dx.
\end{cases}\end{equation}
Here $\phi_j^k(x)$ is the $k$th component of the vector $\phi_j(x)$.  And
$$
u^m(\cdot, 0) = \sum_{i=1}^mg_i^m(0)\phi_i(x),\quad \text{where} \quad
g_i^m(0)=\int_\Omega u_0(x)\phi_i(x)dx.
$$

To finish the proof, we need the following lemma.
\begin{Lemma}\label{Le1}
There exists a weak solution $(\rho^m,u^m,d^m)$ to the problem (\ref{aden})-(\ref{auinit}) in $Q_T=\Omega\times[0,T]$, for any $T\in(0,\infty)$, satisfying
$$
M_1 \leq \rho^m \leq M_2
$$
$$
u^m\in L^\infty(0,T; L^2)\cap L^2(0,T; H^1)
$$
$$
d^m\in L^\infty(0,T;H^1)\cap L^2(0,T;H^2), \quad  |d^m| \leq 1.
$$
Moreover, $(\rho^m,u^m,d^m)$ is smooth in the interior of $Q_T$ and satisfies the energy equality,
\bg\label{energy}
\int_\Omega|\nabla u^m|^2+|\triangle d^m-f(d^m)|^2dx + \frac{d}{dt}\int_\Omega(\frac{1}{2}(\rho^m|u^m|^2+|\nabla d^m|^2)+F(d^m))dx = 0.
\ed
\end{Lemma}

The proof of this lemma is based on an application of the Leray-Schauder fixed point theorem. Let $v^m=\sum_{i=1}^m h^m_i \phi_i \in C^1(0, T; \Hh_m)$. For each $m$ let $\rho^m$ be a solution to
\begin{equation}\label{transport}
\rho + v^m\cdot \nabla\rho = 0
\end{equation}
with initial condition $\rho(\cdot, 0) = \rho_0$. Let $d^m$ be a solution to
$$
d_t + v^m\cdot\nabla d = \Delta d - f(d)
$$
with initial condition $d(\cdot, 0) = d_0$ and boundary condition $d(x, t) = d_0(x)$. The reason the transport equation \eqref{transport} is solvable for $v^m\in C^1(0, T; \Hh_m)$ is due to the regularity of the eigenfunctions of the Stokes operators (cf. \cite{Tem} \cite{La}). Let
$u^m = \sum_{i=1}^m g^m_i \phi_i\in C^1(0, T; \Hh_m)$ be the solution of the system of linear equations
$$
\sum_{i=1}^m\Aa_i^{mj}(t)\frac{d}{dt}g_i^m(t)=-\sum_{i,k}^m\Bb_{ik}^{mj}(t)h_i^m(t)g_k^m(t)
-\sum_{i=1}^m\Cc_i^jg_i^m(t)+\Dd^{mj}(t).
$$
This system of linear equations is solvable because the eigenvalues of matrix of the coefficients $\Aa_i^{mj}(t)$ are bounded from below, since
\begin{equation}\label{inverse}
\Aa_i^{mj}\gamma_i\gamma_j = \int_\Omega \rho |\psi|^2 dx \geq M_1\int_\Omega |\psi|^2 dx\quad \text{where $\psi = \sum_{i=1}^m \gamma_i\phi_i$}.
\end{equation}
Thus we constructed a mapping $\Mm$ with $\Mm(v^m)=u^m$.
The energy estimate (\ref{energy}) play the key to allow one to apply Leray-Schauder fixed point theorem for $\Mm$.\\
\begin{Remark}
We know that from the estimates obtained in Section 2 and 3, we can pass to the limit and conclude that, there exists $(\rho, u,d)$ such that, taking subsequence if necessary,
$$
\rho^m\rightharpoonup \rho \ \mbox { in } L^p(\Omega\times[0,T]) \ \mbox { for any } p\in(1,\infty)
$$
$$
u^m\rightharpoonup u \ \mbox { weakly in } L^2(0,T;H^1) \ \mbox { and strongly in } L^2(0,T;L^2)
$$
$$
d^m\rightharpoonup d \ \mbox { weakly in } L^2(0,T;H^2) \ \mbox { and strongly in } L^2(0,T;H^1).
$$
It follows easily from the above convergences that indeed $(\rho, u, d)$ is a weak solution to the system \eqref{LCD}.
\end{Remark}

\bigskip

\section{Appendix: Frozen Coefficient Method}
For completeness, in this section, we recall an application of the frozen coefficient method by considering the following problem (for more detail see \cite{LS1})
\begin{equation}\label{linear}\begin{split}
\rho v_t-\Delta v+\rho v\cdot\nabla v+\nabla p=f,\\
\nabla\cdot v=0, \ \ v(0)=0, \ \ v|_{\partial\Omega}=0,
\end{split}\end{equation}
in the domain $\Omega\times(0,T)$. Here $\rho$ is a given function in $\Omega\times(0,T)$ satisfying
\bg\label{6den}
0<M_1\leq\rho\leq M_2
\ed
\bg\label{6dend}
|\nabla\rho|\leq M_3, \ \ |\rho_t|\leq M_4
\ed
(the derivatives of $\rho$ may be the generalized ones). And $f$ is a given function in $L^q(\Omega\times(0,T))$ for $q>\sqrt 3$.\\

The frozen coefficient method is now used to prove the following lemma,
\begin{Lemma}
If $v\in W^{1,q}(W^{2,q}(\Omega))$ and $\nabla p\in L^q(\Omega\times(0,T))$ is a solution of problem (\ref{linear}), then for an appropriate constant $a$ depending on $q$
\begin{equation}\label{6Lq}
\|v_t\|_{L^q}+\|\Delta v\|_{L^q}+\|\nabla p\|_{L^q}
\leq C(\|f\|_{L^q}+\|v\|_{L^q}^a),
\end{equation}
where $L^q$ denotes the norm in $L^q(\Omega\times(0,T))$, the constant $C$ depends on $q, \Omega, M_1, M_2, M_3$ and $M_4$.
\end{Lemma}
\proof Let $\zeta_k^\lambda(x)$ be the set of nonnegative functions in $C^2(\R^n)$ depending on the parameter $\lambda\in(0,1)$, forming a partition of unity in $\Omega$:
$$
1=\sum_{k=1}^{N_\lambda}\zeta_k^\lambda(x), \ \ x\in\bar\Omega
$$
with
$$
supp \zeta_k^\lambda(x)\subset\Omega_k^\lambda
$$
where the subdomains $\Omega_k^\lambda$ form a covering of $\Omega$ which satisfy the following properties
\begin{enumerate}
\item the diameters of $\Omega_k^\lambda$ do not exceed $\lambda$;
\item there is a fixed number $l$ such that the multiplicity of the covering of $\Omega$ by the subdomains $\Omega_k^\lambda$ is $\leq l$;
\item $\sum_{k=1}^{N_\lambda}(\zeta_k^\lambda(x))^q\geq \mu>0$, $\forall x\in\Omega$;
\item $|D_x^\beta\zeta_k^\lambda(x)|\leq c\lambda^{-|\beta|}$, $|\beta|=1,2$, $\forall x\in\Omega$.
\end{enumerate}
For simplicity, we consider the case when $\rho$ is independent of $t$. Fix $x_k\in\Omega_k^\lambda$, for each $k=1, 2, ..., N_\lambda$, the vector $v_k(x,t)=v(x,t)\zeta_k^\lambda(x)$ and the function $p_k(x,t)=p(x,t)\zeta_k^\lambda(x)$ are solutions in $\Omega\times(0,T)$ (in the same spaces as $v$ and $p$) to the problem
\begin{equation}\label{klinear}\begin{split}
\rho(x_k)v_{kt}-\Delta v_k+\nabla p_k\\
=f\zeta_k^\lambda(x)+[\rho(x_k)-\rho(x)]v_{kt}-\rho(x)v_k\cdot\nabla v_k+g_k(x,t),\\
\nabla\cdot v_k=v\cdot\nabla\zeta_k^\lambda(x), \ \ v_k(0)=0, \ \ v_k|_{\partial\Omega}=0,
\end{split}\end{equation}
where
\bg\notag
g_k(x,t)=-2\sum_{i=1}^nv_{x_i}\zeta_{kx_i}^\lambda-v\Delta\zeta_k^\lambda+p\nabla\zeta_k^\lambda
+\rho(x)[v\zeta_x^\lambda(v\cdot\nabla\zeta_k^\lambda)+(\zeta_k^\lambda-1)\zeta_k^\lambda v\cdot\nabla v].
\ed
Assuming the right hand side of the first equation in (\ref{klinear}) is the ``free term", making use of the results of Lemma 2.2 in \cite{LS1}, we have the estimate for $v_k, p_k$
\begin{align}\label{vk}
&\|v_{kt}\|_{L^q}+\|\Delta v_{k}\|_{L^q}+\|\nabla p_k\|_{L^q}\\
&\leq C_1(\|f\|_{L^q}+\max_{x\in\Omega_k^\lambda}|\rho(x_k)-\rho(x)|\|v_{kt}\|_{L^q}+\|v_k\cdot\nabla v_k\|_{L^q}+\|g_k\|_{L^q}),\notag
\end{align}
where constant $C_1$ depends on $M_2$.
 Since $\max_{x\in\Omega_k^\lambda}|\rho(x_k)-\rho(x)|\leq M_3\lambda$, and the properties of $\zeta_k^\lambda$, we are able to derive estimate (\ref{6Lq}) from (\ref{vk}), by choosing appropriately $\lambda=(1+M_3)^{-1}\min\left\{1,(2C_1)^{-1}\right\}$ and applying the general Gagliardo-Nirenberg inequality (see \cite{AKM}).
\cbdu
\begin{Remark}As seen in this simple example once we have the bound on the small oscillations the desired estimate follows.
\end{Remark}


{}

\end{document}